\documentclass[5p,10pt]{elsarticle}
 \makeatletter
\def\ps@pprintTitle{%
 \let\@oddhead\@empty
 \let\@evenhead\@empty
 \def\@oddfoot{}%
 \let\@evenfoot\@oddfoot}
\makeatother
\usepackage[utf8]{inputenc}

\usepackage{comment}
\usepackage{graphicx}
\usepackage{amssymb,amsthm,mathtools,amsmath, amsfonts}
\usepackage{hyperref}
\usepackage[nameinlink]{cleveref}

\usepackage{dsfont}
\usepackage{capt-of}
\usepackage{mathtools, cuted}
\usepackage{lipsum, color}
\usepackage{nicefrac}
\usepackage{subcaption}

\usepackage{letltxmacro}

\usepackage{algorithm}
\usepackage[noend]{algpseudocode}

\usepackage[figure]{hypcap}

\usepackage{placeins}
\newenvironment{talign*}
 {\csname align*\endcsname}
 {\endalign}

\usepackage{lineno}
\usepackage{mathtools}
\usepackage{bbm}
\usepackage{float}

\newcommand{\N}{\mathbbm{N}}

\newcommand{\R}{\mathbbm{R}}

\newcommand{\dd}{\, \mathrm d}

\renewcommand{\tt}{\mathbf t}
\newcommand{\uu}{\mathbf u}
\newcommand{\vv}{\mathbf v}
\newcommand{\cc}{\mathbf c}

\newcommand{\cN}{\mathcal N}
\newcommand{\cG}{\mathcal G}
\newcommand{\bdelta}{\boldsymbol{\delta}}
\newcommand{\xn}{x_n}

\newlength\mysinglespace
\newlength\objspace
\newlength\conspace
\newlength\cconspace
\usepackage{todonotes}
\newcommand{\FLil}[1]{\todo[inline,author=FL,color=blue!50,size=\small]{#1}}

\newtheorem{defi}{Definition}

\newtheorem{assm}[defi]{Assumption}

\begin{document}

\begin{frontmatter}


\title{Quality Control in Particle Precipitation via Robust Optimization} 




\cortext[cor1]{Corresponding author (\href{mailto:frauke.liers@fau.de}{\bf frauke.liers@fau.de})}

\author[label0]{Martina Kuchlbauer}
\author[label0]{Jana Dienstbier}
\author[label3,label4]{Adeel Muneer}
\author[label0]{Hanna Hedges}
\author[label3]{Michael Stingl}
\author[label0]{Frauke Liers\corref{cor1}}
\author[label3,label4]{Lukas Pflug}

\address[label0]{Friedrich-Alexander-Universität Erlangen-Nürnberg
  (FAU), Optimization under Uncertainty and Data Analysis, Departement of Data Science, Cauerstraße 11, 91058 Erlangen, Germany}
\address[label3]{Friedrich-Alexander-Universität Erlangen-Nürnberg (FAU), Chair of Applied Mathematics (Continuous Optimization), Cauerstraße 11, 91058 Erlangen, Germany}
\address[label4]{Friedrich-Alexander-Universität Erlangen-Nürnberg (FAU), Central Institute for Scientific Computing, Martensstr. 5a, 91058 Erlangen, Germany}

\begin{abstract}
  We propose a robust optimization approach to mitigate the impact of uncertainties in particle precipitation. Our model of particle synthesis incorporates, as partial differential equations, nonlinear and nonlocal population balance equations. 
  The optimization goal is to design products with desired size distributions. Recognizing the impact of uncertainties, we extend the model to robustly hedge against them to ensure  tailored particle sizes. For the robust problem, we enhance an adaptive bundle framework for nonlinear robust optimization integrating the exact method of moments approach for the population balance equations. Computational experiments  focus on uncertainties in precursor inflow rates, which greatly influence the resulting product's quality. Using realistic parameter values for quantum dot synthesis, we demonstrate the algorithm's efficiency and find that unprotected processes fail to achieve desired particle sizes, even for small uncertainties, highlighting the need for  robust processes. The latter  outperforms the unprotected process concerning the product's quality in perturbed scenarios.
\end{abstract}



\end{frontmatter}

\section{Introduction}
Designing (nano)particles with optimal properties involves numerous challenging tasks, such as establishing suitable experimental setups and processes to produce the desired  products. 
Indeed, the optimization of particle size and shape has a wide
variety of applications, ranging from
development of cosmetics and pharmaceutical usage for oral delivery of
peptide drugs to the performance improvement of lithium batteries
\cite{naito2018nanoparticle}. \\
Achieving optimal results in particle synthesis is particularly challenging due to the presence of uncertainties in, e.g., material parameters, temperature, particle nucleation, and growth, to name just a few, that can have a detrimental impact on the quality of the final product. 
Since even small disturbances can lead to a deterioration of product quality
\cite{mastronardi2012size}, it is essential to ensure that possible uncertainties or fluctuations affect the product as little as
possible. Neglecting uncertainties in particle design can have severe consequences, as even slight variations in their dispersity can significantly impact the performance of the final product. 

As a consequence, 
robust mathematical optimization offers an approach to overcome this issue by enabling the development of optimized processes that are best possible, even under uncertain conditions.
Developing optimization frameworks for particle synthesis and employing suitable methods is of utmost importance in achieving efficient and practical manufacturing processes. Robust design is particularly crucial in industrial contexts to guarantee desired characteristics such as absorption, emission properties, quantum yield, and maximum mass while optimizing production costs. Uncertainties, encompassing factors like solubility, flow rates, and temperature, profoundly influence nanoparticle synthesis even at the laboratory level. These uncertainties, coupled with challenges in characterization, significantly affect the overall quality of the product. Consequently, effectively managing and mitigating uncertainties and characterization issues becomes imperative to attain nanoparticles of high quality. For hedging against uncertainty,  the currently very active research area of robust optimization is very well suited. 
However, robust optimization adds additional challenges to the optimization of process designs. In this work, we contribute to this area by introducing a practically efficient robust solution algorithm that is able to determine optimal processes for particle precipitation. We show its effectiveness for robust processes. 

\subsection*{Our Contribution}
We present a quality control system for particle precipitation modeled by population balance equations (PBEs). For this nonlinear, nonlocal class of partial differential equations (PDEs), no analytical closed-form solutions are available. However, semi-analytical solutions can be used by the exact method of moment (eMoM) approach. For hedging against uncertainties, we integrate eMoM into a novel adaptive bundle framework for robust (nonlinear) optimization. We use it in particular for optimizing particle size distributions (PSDs) under uncertain precursor inflow rates as the latter has a strong impact on the result of the synthesis process. 
The computational results convincingly show the need for robust protection. In addition, it turns out that a robust process yields products of desired sizes while keeping the cost of robust protection very low. 

\subsection*{Related Literature} 

The impact of dispersity on product performance is particularly critical for small nanoparticles and new materials, where slight variations can significantly affect properties such as band gap and quantum yield \cite{alivisatos1996semiconductor,talapin2001highly,micic1994synthesis,burda2005chemistry,viswanatha2004understanding,mastronardi2012size}. However, the mathematical optimization of robust nanoparticle design is a yet largely unexplored area, and we aim at contributing to filling this gap. Mathematical optimization provides a unique opportunity to establish optimized processes that ensure product quality under uncertainty which is, in particular, relevant for up-scaling the processes.
The PBEs that are used to model particle precipitation are nonlinear PDEs for which typical solution approaches are finite-volume type methods \cite{qamar2009solution, gunawan2004high},
methods based on characteristics \cite{fevotte2010method, ur2014application}, and
methods based on moments of the solution such as MoM \cite{schwarzer2006predictive}, DQ-MoM \cite{marchisio2005solution}, and eMoM \cite{pflug2020emom}
which is based on theoretical insights from 
\cite{keimer2017existence,keimer2018nonlocal}.
In this work, we set up an optimization problem for particle precipitation that uses the eMoM representation. We also go one step further and hedge the optimized process against uncertainties.

Indeed, in many applications, uncertainties largely impact the quality of the obtained solutions and can easily render an obtained result useless. Typical reasons are that uncertainties can lead to high costs, to a low-quality solution, or the solutions obtained may even be infeasible as they do not satisfy the side constraints that are affected by the uncertainty. Several solution approaches have been presented in the literature. Stochastic optimization is a well-known paradigm used to address optimization problems with uncertain parameters. This approach involves the use of random variables, such as random objective functions, random constraints, and stochastic quantities like expected values or higher moments. In various literature, such as \cite{prekopa1998SO}, \cite{birge2006introduction}, and \cite{Shapiro2003}, these concepts are discussed in detail. It should be noted that stochastic optimization requires knowledge of the distributions of uncertain parameters, which is often not available. In addition, protection against uncertainties is only guaranteed in a probabilistic sense and with a certain probability. In contrast, in the case of nanoparticle design, it is essential to guarantee quality requirements with certainty. Moreover, the underlying distributions are themselves uncertain. This requires a different and robust approach. 

Robust optimization does not require knowledge of the underlying probability distributions. It addresses these needs by predefining so-called uncertainty sets against which protection is sought. The resulting robust problem ensures feasibility for all realizations of the uncertain parameters within the uncertainty set and optimizes the guaranteed cost. For continuous uncertainty sets, robust optimization models are typically represented as semi-infinite problems (SIPs) with finitely many variables and infinitely many constraints. In order to obtain algorithmically tractable robust counterparts of these SIPs, elegant theoretical as well as algorithmic concepts have been developed. For a focus on robust linear and convex optimization where - among others - reformulations to tractable counterparts have been presented, we refer to the books \cite{Ben-Tal:RobustOptimization2009}, \cite{denhertog2022}. 
Robust optimization approaches have been developed and applied in the context of process optimization, e.g., in 
\cite{zukui2008_robust_process, wiebeMisener2019_robopt_pooling, koller2018_stoch_back-off}. 
In particular, for a restricted synthesis process assuming seeded growth and linearized PBEs, a convex robust optimization task could be established, together with an algorithmically tractable reformulated robust counterpart in \cite{dienstbier2022robust}. 

More generally, when including nucleation and growth and using the general nonlinear PBE, we face nonlinear and non-convex robust optimization problems. For this challenging class of optimization problems, reformulation approaches are typically not available, and
approximation schemes, such as \cite{Diehl_ApproxTech, Houska2013nonlinear}, or decomposition methods are a potential choice. Indeed, only a few general approaches are known for this setting. In this work, we enhance a specific decomposition algorithm. The latter consists in a novel adaptive bundle method that has been developed as a general framework for nonlinear robust optimization in \cite{kuchlbauer2022adaptive}. The method has been integrated into an outer approximation scheme for discrete-continuous decisions in \cite{kuchlbauerOuter}. In this work, we make it concrete for optimal particle precipitation. To this end, we present the adaptive bundle method, which integrates the eMoM approach and uses it for a practically efficient solution of the underlying PBE.  


\subsection*{Outline} 
We introduce the setting in Section \ref{sec:nominaloptimization} and summarize the mathematical model for particle synthesis via PBEs including nucleation and growth. We revisit the solution approach via eMoM \cite{pflug2020emom} and also introduce the nominal optimization problem for particle precipitation. 
Section \ref{sec:robustification} extends the model by considering mass uncertainties and presents the robust optimization problem that hedges the precipitation process against uncertainties. Section \ref{sec:bundle} summarizes the adaptive bundle method \cite{kuchlbauer2022adaptive} and presents how to integrate the eMoM that is used for the solution of the PBE. Computational results for the unprotected as well as for the robust synthesis process are presented and discussed in \ref{sec:results}. We conclude with a discussion and an outlook for future research.

\section{Optimized Synthesis modeled by Population Balance Equations}\label{sec:nominaloptimization}

The evolution of the distribution of particle sizes is based on \textit{nucleation} and
\textit{growth}, which are driven by the supersaturation of the reactant 
concentrations. Particle precipitation can be described by a PBE. With the diameter of crystals with spherical shape denoted by $x$, process time by $t$ and the particle size distribution denoted by $q$ is a solution of this PBE, see \cite{pflug2020emom} and \cite{Gilch2021}. 


An important goal in particle design is to establish the best possible synthesis processes with respect to an objective functional, such as maximizing the amount of product that satisfies specific size or, more generally, quality criteria. We will later use 
the time-dependent inflow rate of precursor as the control variable because it has a significant effect on the final particle size distribution. 


\subsection{Population Balance Equation}

The precursor concentration $c_{pre}$ and the total concentration of reduced educt species $c_{tot}$ are governed by the control $v$, i.e. the precursor inflow rate, and prescribed by the following set of ordinary differential equations: 
\begin{equation}
\begin{aligned}
    \partial_t c_{pre}(t) &= v(t) - k_r c_{pre}(t),      & c_{pre}(0) &= 0,\\
\partial_t c_{tot}(t) &= k_r c_{pre}(t), &c_{tot}(0) &= 0,
\end{aligned}
\end{equation}
where $k_r \in \R_{>0}$ is the reduction rate of the precursor species.

Let $\mathcal{G}_{0}$ be the concentration-dependent and $\mathcal{G}_{1}$ the size-dependent part of the growth rate, $q_{0}$ the initial number density distribution, $\mathcal{N}$ the concentration-dependent nucleation rate, and $\xn \in \mathbb{R}_{>0}$  the nucleation size of particles. 
The solution of the following system of equations then gives the number density distribution $q$:
\begin{equation}
\label{eq:PBE}
\begin{aligned}
\partial_t q(t,x)  + \partial_x\big(\mathcal{G}_{0}( c(t) )\mathcal{G}_{1}(x)q(t,x) \big) &= 0, \\
q(0,x) & = q_{0}(x), \\
\mathcal{G}_{0}\left( c(t) \right)\mathcal{G}_{1}(\xn)q(t,\xn) &= \mathcal{N}\left(c(t) \right).\\
\end{aligned}
\end{equation}
We consider only nucleation and no preexisting particles in this contribution, thus $q_{0} \equiv 0$. The total concentration $c_{tot}$, the solution in the liquid phase $c$, and the PSD $q$ are connected by the mass balance as follows: 
\begin{align}
    c(t) = c_{tot}(t) -  \tfrac{\rho\pi}{6V}\int_{\xn}^{\infty} x^{3}q(t,x)\, \mathrm dx.
\end{align}
where $V$ denotes the
total volume of the system, $\xn$ represents the nucleation size, while $\rho > 0$ signifies the nanoparticle density. 
As the system of equations governing the evolution of the precursor concentration and the total educt concentration can be solved analytically, we obtain:
\begin{align}
    c(t) = \int_{0}^t \!\!\big(1-e^{-k_r(t-\tau)}\big) v(\tau) \, \mathrm d\tau  -  \tfrac{\rho\pi}{6V}\!\!\int_{\xn}^{\infty} \!\!\!\! q(t,x)x^{3}\dd x.
\end{align}
\noindent
In addition, it is assumed that the size-dependent growth rate
$\mathcal{G}_{1}$ lies in the class of growth kinetics of the form
$\mathcal{G}_{1} = x^{\beta}$, for $\beta \in \R \setminus \{1\}$. The examples presented in this paper focus on diffusion-limited growth, i.e. $\beta = -1$. However, we would like to point out that the generality of the presented approach also covers other relevant growth kinetics in nanoparticle precipitation, such as those shown in \cite{thanh2014mechanisms}.

\subsection{Exact Method of Moments}
The exact method of moments, as derived in \cite{pflug2020emom}, approximates the evolution of the concentration in the liquid phase $c$, instead of solving the PBE. For $t\in \left[0,t_{n}\right]$, we can define:
\begin{equation}
\label{defi:vol}
\begin{aligned}
c(t) &= \int_{0}^t \big(1-e^{-k_r(t-\tau)}\big) v(\tau) \, \mathrm d\tau \\
&\qquad - \tfrac{\rho\pi}{6 V}     \int_{0}^{t} \mathcal{N}\left( c(\tau )\right) \psi \left[\tau ,\xn\right](t)^{3}\, \mathrm d\tau,
\end{aligned}
\end{equation}
where
\begin{equation}
\label{psim}
\begin{aligned}
\psi \left[\tau ,x \right](t) &:= \left(x^{1-\beta} + (1-\beta)\int_{\tau}^{t}\mathcal{G}_{0}(c(s))\, \mathrm ds\right)^{\frac{1}{1-\beta}}. \\
\end{aligned}
\end{equation}
Here, $\psi \left[\tau,x \right](t)$ describes the size of a particle at the time $t$, with the condition that this particle was of size $x$ at the time $\tau$. As shown in \cite{pflug2020emom}, given a solution $c$ of \eqref{defi:vol}, we obtain at a solution formula for the PBE as follows. 

For $(t,x) \in \R_{>0}\times\mathbb{R}_{\geq \xn}$, the solution of the PBE in \cref{eq:PBE} is given by:
\begin{equation}
\label{eq:qm}
\begin{aligned}
q(t,x) &= x^{-\beta} \begin{cases}
\frac{\mathcal{N}\left( c(\tau_{t,x}) \right)}{\mathcal{G}_{0}\left(  c(\tau_{t,x}) \right)} & \text{ for } x\leq \psi\left[0 ,\xn\right](t) \\
 0 & \text{otherwise} \end{cases},
\end{aligned}
\end{equation}
with $c$ satisfying equation \eqref{defi:vol} and 
\begin{equation}
\tau_{t,x} := \max\left\{\nu\in\left[0,t\right] : \psi\left[\nu,\xn\right](t) = x\right\},
\end{equation}
for all $t\in\left[0,t_{n}\right]$ and $x\leq \psi_{m}\left[0,\xn\right](t)$.
\noindent

In order to formulate desired properties of the final PSD in terms of their moments, e.g. the mean and the variance of the final PSD, the moments of the PSD can also be given directly from the concentration profile. This means that the p-th moment $m_{p}$ of the number density distribution $q$ at time $t \in \R_{\geq 0}$ is given by:
\begin{equation}
\label{eq:moments}
\begin{aligned}
m_{p}(t) :&= \int_{\xn}^{\infty}q(t,x)x^{p}\, \mathrm dx, \\
&= \int_{0}^{t} \mathcal{N}\left(  c(\tau )\right) \psi \left[\tau ,\xn\right](t)^{p}\, \mathrm d\tau, 
\end{aligned}
\end{equation}
where $\psi$ is again given by \eqref{psim}.



\subsection{An Optimization Model for Particle Precipitation}
In order to establish an optimization problem for synthesis processes leading to tailored PSDs, we aim here to optimize the number-weighted mean and variance of $q$ at the end of the synthesis process. More specifically, the goal is to control the precursor feed to the system so that the deviation of the mean of the distribution from a desired mean and the deviation of the standard deviation of the distribution from a desired standard deviation are minimized.  

Let $\bar{m}_{0}$ be the desired number-weighted mean and assume that we want to control the standard deviation to zero. Then, the objective functional $J$ is given by:
\begin{equation}
\label{eqn:COST}
J := \omega_1 \left(\tfrac{m_1(T)}{m_0(T)}- \bar m_0\right)^2 \!\! + \omega_2\left(\tfrac{m_2(T)}{m_0(T)} - \left(\tfrac{m_1(T)}{m_0(T)}\right)^2\right).
\end{equation}
Note that the squared standard deviation, i.e. the variance, of a random variable $X$ is given by $E[X^2]-E[X]^2$ while the weights $\omega_1,\omega_2 \in \R_{\geq 0}$ scalarize the two criteria. 

\subsubsection{Discretization}\label{sec:discr}
In order to numerically approximate the solution of the PBE and its moments, we use a discretization analoguously to \cite{pflug2020emom} for \eqref{defi:vol} and \eqref{eq:moments}. For this, we discretize the time horizon into $N_t \in \N_{>0}$ time-points $0 = \tt_0 < \tt_1 < \ldots < \tt_{N_t} = T$ and define the length of each time-interval by $\bdelta_k := \tt_{k+1}-\tt_{k}$. Based on this, we approximate the control $v$ by its values at these time-points, i.e. $\vv_i = v(\tt_i) \ \forall i$. In this way, we can approximate the total mass in the system $c_{tot}$ as follows:
\begin{align}
\cc_{tot,k}&=\sum_{\ell=1}^{k-1} \vv_\ell \Big(\bdelta_\ell+\tfrac{1-e^{k_r\bdelta_{\ell}}}{k_r e^{k_r (\tt_k-\tt_\ell)}}\Big) \label{eqn:DIS1},
\end{align}
and based on this the evolution of the approximated concentrations $\cc_k \approx c(\tt_k)$ as follows:
\begin{equation}
\begin{aligned}
\cc_{k+1} 
&= \cc_k +\gamma_1 \big(\cc_{tot,k+1}-\cc_{tot,k}\big)  -\gamma_2 \cN(\cc_{k})\bdelta_k\cG_{k,k}^3\\
& \qquad -\gamma_2 \sum_{\ell=1}^{k-1} \cN(\cc_{\ell}) \bdelta_\ell \big(\cG_{\ell,k}^3 -\cG_{\ell,k-1}^3\big) \label{eqn:DIS2}, \\
\cc_1 &=0,
\end{aligned}
\end{equation}
with the approximated characteristics:
\begin{align}
\cG_{\ell,k} := \Big(\xn^{1-\beta}+(1-\beta) \sum_{m=\ell}^{k} \cG_0(\cc_m)\delta_m\Big)^{\frac{1}{1-\beta}}.
\end{align} 
All of these discretizations are based on dividing the integral involved into the time intervals and then approximating the integrand by its value at the lower bound of each interval.

Combining the derived governing equations of the discretized concentration $\cc$ based on the discretized control $\vv$, we obtain:
\begin{align}\label{eqn:DISCRSTATEOP}
    F(\vv,\cc) = 0,
\end{align}
with $F : \R^{N_t} \times \R^{N_t} \mapsto \R^{N_t}$. Let $C : \R^{N_t} \mapsto \R^{N_t}$ be the mapping from the control $\vv$ to the respective concentration $\cc$, i.e. $F(\vv,C(\vv))= 0$.

\subsubsection{Nominal Optimization Problem}\label{sub:nominal}
Let \begin{align}\label{eqn:ADMSET}
\mathcal{V} := \left\{\vv \in [l,u]^{N_t} : \textstyle\sum\limits_{k = 1}^{N_t} \bdelta_k \vv_k = V_{tot}\right\},
\end{align}
be the set of admissible controls, with lower and upper bounds $0\leq l< u$ and total added precursor $V_{tot}$. The optimization problem can be formulated in the following way:
\begin{equation}\label{nominalproblem}
\min_{\vv\in \mathcal{V}}J(\vv,C(\vv)),
\end{equation}
where, as stated above, $C(\vv) \in \mathbb{R}^{n}$ satisfies
$F(\vv,C(\vv)) = 0$.

To solve this optimization problem by first-order methods, we apply the implicit function theorem to $F$ to obtain the derivative:
$$\nabla C(\vv) = -\nabla_2 F(\vv,C(\vv))^{-1}\nabla_1 F(\vv,C(\vv)),$$ 
and consequently
\begin{align}\label{eqn:GRAD}
&\nabla_\vv J(\vv,C(\vv)) = \nonumber\\
&\qquad \nabla_1 J(\vv,C(\vv)) + \nabla_2 J(\vv,C(\vv))\nabla C(\vv)).
\end{align}

\section{Particle Precipitation under Uncertainties\label{sec:robustification}}
We acknowledge that in a synthesis process set up in a laboratory, it is not possible to exactly control precursor inflow rates. This means that the resulting uncertainties can have a significant impact on the outcome of the process. Therefore, our focus will be primarily on adjusting this parameter. As a consequence, we will formulate a robust optimization problem. To solve it, we will extend an adaptive bundle framework from \cite{kuchlbauer2022adaptive} and combine it with the exact methods of moments.
To model uncertainties in the mass flow, we use a multiplicative uncertainty set that captures relative parameter variations. The latter allows us to relate the parameter variations to their magnitude, taking into account that a smaller parameter may have smaller variations compared to a larger parameter. 
We define our model according to Section \ref{sec:nominaloptimization} in a time-discretized setting.
The aim is to find a solution that is robustly protected against relevant deviations from a nominal, i.e. typical, scenario. 


Thus, the uncertainty set naturally includes the nominal scenario (without disturbance), which corresponds to a factor of $1$ for each time point. A commonly used uncertainty set is 
given by a $\Vert\cdot\Vert_\infty$-bounded function of the time interval, resulting in a box-shaped uncertainty set. Assume that the worst cases scenario if of "bang-bang" type, i.e., the worst cases scenarios can be found in a vertex of the uncertainty set. If we further assume that pushing the inflow rate in the nucleation phase at the beginning of the nanoparticle precipitation and the growth phase in the second part of the synthesis to opposite boundaries causes the strongest changes in the objective functional we end up with a monotone uncertainty with values of either the upper or the lower bound.  To also encounter for the nominal case, i.e. uncertainty $\uu = 1$, we consider all monotone perturbations with function values $u_l,1,u_u$ and a single jump within the process time. 
\begin{figure}[h]
	\includegraphics[width=0.5\textwidth]{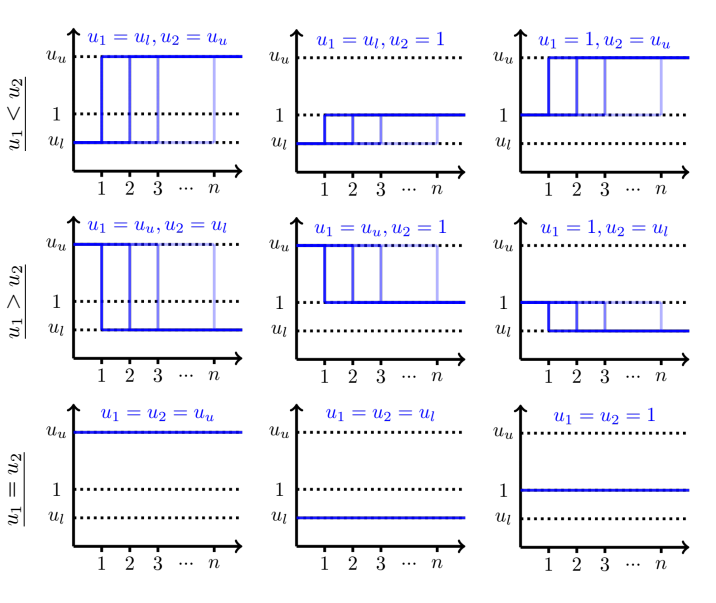}
	\caption{Conceptual visualization of the used uncertainty set (\ref{eq:U})}
	\label{fig:u_set}
\end{figure} 
In more detail, we define the uncertainty set as follows. 

Let $n$ be the number of discretization steps and $u_l , u_u \in \mathbb{R}$ such that $0 < u_l < 1 < u_u < \infty$. The uncertainty set is given by:
\begin{align}
\label{eq:U}
\mathcal{U}:= \Big\{ &\uu \in \mathbb{R}^{n} | u_1, u_2 \in \{u_l, 1, u_u\}, j \in \left\{ 1,...,n \right\}, \\
& \uu_i = u_1 \text{ for } i\leq j, \uu_i = u_2 \text{ for } i> j \Big\}.
\end{align}
This uncertainty set then represents the basic building block for more general uncertainty sets where more than one change in uncertainty is manifested over the full-time period. The uncertainty considered here can easily be extended in the future to account for multiple changes in the factor over time.

Next, we obtain a robust optimization model as follows.
Let $\mathcal{V}$ be the compact set of admissible discretized controls as defined in \cref{eqn:ADMSET} and the cost function $J$ be chosen as in \cref{eqn:COST} depending on the discretized state operator $F$ and the discretized concentration function $C$ introduced in \Cref{sec:discr} above.
Then the robust optimization problem is:
\begin{equation}\label{eq:rob_opt_syn}
\begin{aligned}
\min_{\vv\in \mathcal{V}} \max_{\uu \in \mathcal{U}} \; & J(\uu\bullet \vv,C(\uu\bullet \vv)),
\end{aligned}
\end{equation}
where $\uu\bullet \vv$ denotes the entrywise product of $\uu$ and $\vv$, i.e.,
$$
(\uu\bullet \vv)_i = \uu_i \vv_i \qquad \forall i=1,\ldots,n.
$$
Moreover, as before, the concentration $C(\uu\bullet \vv)$ is implicitly defined as the unique solution of the discretized state equation:
\begin{align}
    F( \vv \bullet \uu , C(\uu\bullet \vv)) = 0.
\end{align}
Since this problem is nonlinear and non-convex, there is no standard procedure for a practical or theoretical solution. In the next section we will focus on its solution via an adaptive bundle framework by \cite{kuchlbauer2022adaptive}.

\section{Robust Optimization via Bundle Method for Nonconvex Optimization \label{sec:bundle}
}


We tackle the robust PSD problem \eqref{eq:rob_opt_syn} by the adaptive bundle method from \cite{kuchlbauer2022adaptive}. This bundle method is applicable to robust optimization problems of the form:
\begin{equation}
\label{eq:minmax}
\min_{\vv\in\mathcal{V}}\max_{\uu\in \mathcal{U}} h(\vv,\uu),
\end{equation}
with a compact uncertainty set $\mathcal{U}$ and a cost function $h:\mathbb{R}^{n+m} \rightarrow \mathbb{R}$ that is  locally Lipschitz continuous and lower $C^1$. 
Choosing $h(\vv,\uu) := J(\uu\bullet \vv,C(\uu\bullet \vv))$, it is readily seen that the robust PSD \eqref{eq:rob_opt_syn} fits into this setting. 

More generally, the adaptive bundle method solves the problem \eqref{eq:minmax} by treating it as a nonsmooth minimization problem with the optimal value function: $$\hat{h}(\vv) := \max_{\uu\in \mathcal{U}} h(\vv,\uu),$$ as an objective function. 
An evaluation of the objective function, therefore, requires solving the inbuilt maximization problem. 
This maximization problem aims to find the worst-case realization of the uncertain parameter. 
We thus call it the adversarial problem. For the PSD problem, the adversarial problem is to find the value of the uncertain parameter $\uu$ that maximizes the deviation from the desired  volume-weighed mean and standard deviation as defined in \cref{eqn:COST}. 

The bundle method requires access to approximate solutions to the adversarial maximization problem in the following sense.
\begin{assm}
We assume that for every $\vv\in \mathcal{V}$ and every choice of $\epsilon_{f} > 0$, we have access to a $\uu_\vv \in \mathcal{U}$ such that $h(\vv,\uu_\vv) \geq \max_{\uu\in \mathcal{U}}h(\vv,\uu) - \epsilon_f$, to the corresponding function value $h(\vv,\uu_\vv)$ and to a subgradient $g_\vv \in \partial_\vv h(\vv,\uu_\vv)$.
\end{assm}
Therefore, it is necessary to compute an approximate worst-case for every decision and tolerance. The bundle method does not require any assumptions on the structure of the underlying optimization models, which allows it to be used for nonlinear robust optimization. In addition, it is also not restricted to finite-dimensional optimization problems but can be extended more generally to robust optimization of models containing PDEs, as we do here for robust particle precipitation, including PBE. 
According to the assumption, the method requires an approximate worst-case realization for a given tolerance $\epsilon_f$:
\begin{equation}
\uu_\vv \in \left\{ \uu\in \mathcal{U} \ | \ h(\uu,\uu_\vv ) \geq \max_{\uu\in \mathcal{U}} h(\uu,\uu_\vv ) - \epsilon_f \right\}.
\end{equation}
For the first-order information required by the bundle method, we can use an approximate subgradient:
\begin{equation}
g_\vv \in \partial_\vv h(\vv,\uu_\vv).
\end{equation}
Note that in this paper the error $\epsilon_f$ can be chosen to be $0$ since we only allow jumps of the uncertainty function at a finite number of points, which are given by the time discretization. As a consequence of this, the required subgradient is computed directly by evaluating \cref{eqn:GRAD} in the perturbed control $(\vv,\uu_\vv)$. 

We mention, however, that the adaptivity in the bundle approach can be exploited in future work to speed up computations for very fine time discretizations, or, even more interestingly, if a continuous uncertainty set were used. An example of the latter would be obtained, if we allowed our monotonous uncertainty functions to jump at arbitrary points. As will be shown in the computational results, the approach used for our purposes already leads to interesting results, so these algorithmic improvements are postponed to future research. 

We will now explain the algorithmic concept of the bundle method. The iterations are divided into inner and outer loops. The outer loops generate serious iterates while the inner loops generate trial iterates around the current serious iterate. These are then potential candidates for the next serious iterate. To distinguish between the inner and outer iterations in the algorithm and in the following statements, we denote the inner loop counter by $k$ and the respective trial iterations by $\vv^k$, as well as the outer loop counter by $j$ and the respective  serious iterations by $\vv_j$. The algorithm also contains the proximity control parameter $\tau_k$, which depends on the current quality of the working model. Details can be found in \cite{kuchlbauer2022adaptive}. We initiate the algorithm by running the nominal optimization. The solution is our first serious iteration $\vv_1$.

New trial iterates are generated by optimizing a convex working model. The working model around a serious iterate $\vv$ is constructed as a maximum of cutting planes, i.e. as a piecewise linear model. To do this, the method requires approximate function values and subgradients. For each trial iteration, we therefore solve the adversarial maximization problem up to a given error. The required  error bound is an adaptive one and depends on the current distance between iterations. 
\begin{algorithm}[h]
\caption{Adaptive bundle method \cite{kuchlbauer2022adaptive}} \label{alg:bundle}
\begin{algorithmic}[1]

    \State Initial  $\vv_0$: Run nominal optimization. Set $j = 1$
    \While{dist$(0 , \tilde{\partial}_{a}\hat{h}(\vv_j ) > \tilde{\epsilon}$} 
    	\State Start inner loop with serious iterate $\vv_j$; set $k = 1$.
        \State Find new trial iterate $\vv^k$ by solving  
        \begin{equation*}
        \min_{\vv^k \in \mathcal{V}} \Phi_k (\vv^k , \vv_j ) + \frac{\tau_k}{2}|| \vv^k - \vv_j ||^2
        \end{equation*}
        
        \State Approximate $\hat{h}(\vv^k)$ 
        \If{ accept} 
        \State $\vv_{j+1} \leftarrow \vv^k$
        \Else
        \State Generate cutting plane with $g_k \in \tilde{\partial}_a \hat{h}(\vv^k )$
        
        \State Build a new working model 
        \State Increase $k$ and go to Step 5.
        \EndIf

    \State Increase j.
\EndWhile

\end{algorithmic}
\end{algorithm}

\section{Implementation\label{sec:implementation}}
The adaptive bundle method was implemented in MATLAB. The experiments were performed on a machine with an Intel Core i7-8550U (4 cores) and 16GB of RAM. For the computation of our robust solutions, we encountered running times between $5$ and $20$ minutes.

The initial value for the bundle method was chosen as the result of the nominal optimization. 
The algorithm used for this is an adjusted version of \cite{Gilch2021}. \\
In step 5 of the bundle algorithm, we solve the adversarial maximization problem over the (discretized) uncertainty set by a straight-forward enumeration of all possible discretized values. 

\section{Computational Results}\label{sec:results}

In this section, we present the numerical results for optimizing the particle synthesis process. 

Table 1 shows the parameter values that were used in the optimization model. To choose realistic values, we assume that the quantum dots (QDs) should have a mean value of 4 nm, and that the process should end after 10 minutes. The other values can be found in the table. The optimization goal is to minimize the quadratic difference of the mean size and the variance from a given value, where their effects are weighted by two parameters $\omega_1, \omega_2$, see (\ref{eqn:COST}). 

In the computational evaluations presented here, we weight both contributions to the objective equally, see table.
However, it turns out that the results are qualitatively similar even for different choices, as can be seen in Appendix \ref{sec:appendix-computations}, where we present results for a different choice.  

\begin{table}[h]
\begin{center}
\begin{tabular}{ l | c | c | c } 
 \textbf{name} & \textbf{symbol} & \textbf{value} & \textbf{unit} \\
 \hline
 desired NP size & $\bar{y}$ & 4 & nm\\ 
 final control time & $t_n$ & 10 & min\\ 
 nuclei size & $x_{\star}$ & 1 & nm\\
\# control points & n & 100 &  -\\
growth law exponent & $\beta$ & -1 & -\\
nanoparticle density & $\rho$ & 1 & -\\
obj. func. factor size &$\omega_1$ & 1 & -\\
obj. func. factor var. &$\omega_2$ & 1 & -\\ 
\end{tabular}
\label{tab:optim-parameters}
\end{center}

\caption{Optimization parameters as used in the computational results.} 
\end{table}

In particular, we focus on two optimized solutions: the \textit{nominal} optimum synthesis process (\ref{nominalproblem}) from Section \ref{sub:nominal} without protection against uncertainties, and the robust process (\ref{eq:rob_opt_syn}) which is designed to hedge against uncertainties. We are interested in evaluating these processes when used in the \textit{nominal}, i.e., unperturbed, scenario which represents ideal laboratory conditions where all $\uu_i$ values are set to 1. In addition, we study perturbed scenarios where uncertainties come into play. To study extreme cases, we analyze the worst scenario for both the nominal and robust synthesis processes. The \textit{nominal worst case} corresponds to the uncertainty realization that results in the largest objective function value for the nominal process. Similarly, the \textit{robust worst case} represents the uncertainty where the robust process achieves the worst objective function value.

\begin{figure}[h]
\centering
\includegraphics[width=0.49\textwidth]{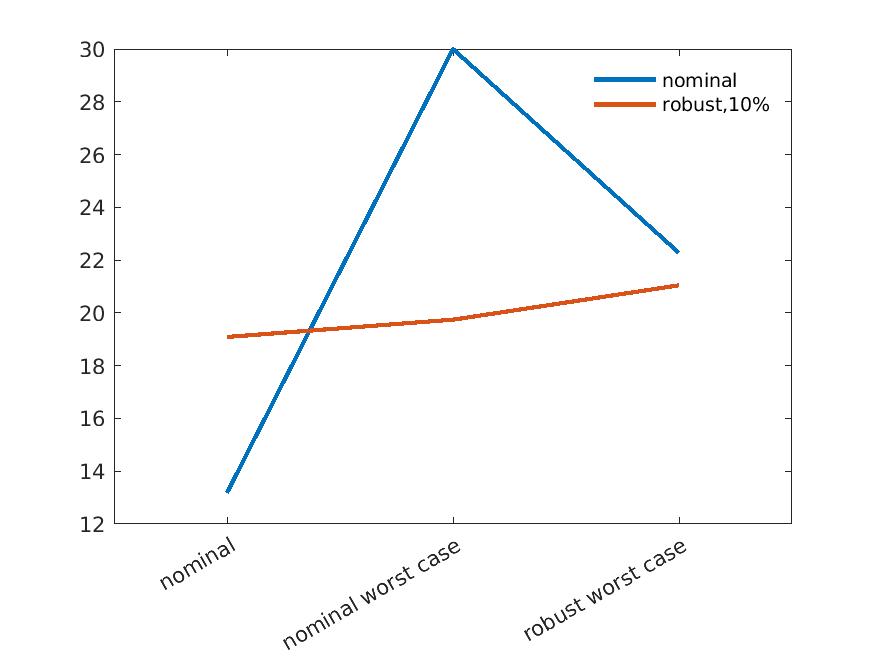}
\caption{
Objective function values for the uncertainty set with $[u_l, u_u] = [0.9, 1.1]$. Lines are guides for the eyes only.}
\label{fig:rob1}
\end{figure} 
Figure \ref{fig:rob1} illustrates the objective function values for the nominal and robust processes, where the robust process is obtained with a maximum uncertainty of 10\% ($[u_l, u_u] = [0.9, 1.1]$). The results are presented for three scenarios: the nominal, the nominal worst case, and the robust worst case.

In the nominal scenario, as expected, the nominal process, specifically optimized for this case, outperforms the robust one. However, the situation changes significantly when considering perturbed scenarios, particularly when the uncertainty manifests within the 10\% uncertainty set to its nominal worst case. Assuming that this scenario manifests itself, the nominal objective value then deteriorates by more than a factor of two, when compared to the nominal scenario. This indicates that the nominal process is highly sensitive to uncertainty. When compared to the robust process, it is still approximately $\frac{1}{3}$ worse than the robust one. Conversely, the robust process demonstrates stable objective function values that exhibit only slight changes from the nominal to the robust worst-case scenario. 
Also in the robust worst-case scenario, the robust process performs slightly better than the nominal one.   

\begin{figure}[h]
\centering
\includegraphics[width=0.49\textwidth]{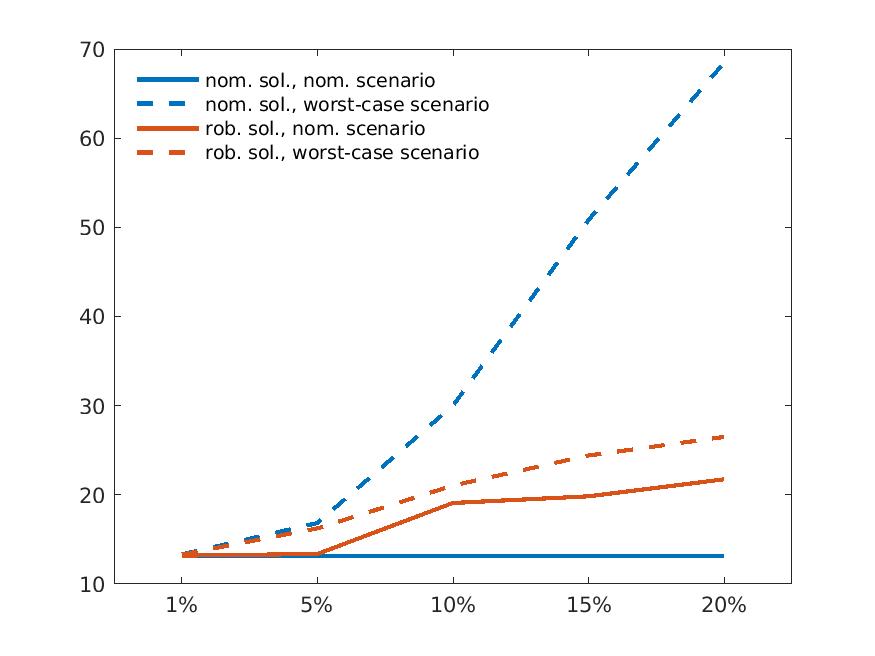}
\caption{
Objective function values for both the nominal and the robust process, are evaluated in different scenarios, as a function of the size of the uncertainty set. The size of the uncertainty is given in percent deviation from the nominal scenario.}
\label{fig:rob2}
\end{figure} 
Figure \ref{fig:rob2} displays the objectives for the nominal process, evaluated both for the nominal and the nominal worst-case scenario, and for the robust process, evaluated both for the nominal and the robust worst-case scenario. It is clear that the nominal optimum process (blue line) stays at a constant objective value in the nominal, i.e., unperturbed, scenario, irrespective of the size of the uncertainty set. However, when the nominal process is evaluated in the worst nominal scenarios, the objective value of the latter significantly increases to roughly a factor of 5 for the largest uncertainty of 20\%. This means that if this uncertainty scenario is realised, the quality of the synthesised product obtained will be significantly reduced, as the mean and variance of the products obtained will be significantly different from the given values. Comparing these results with the robust process (red lines), the obtained robust objectives are significantly better than those of the nominal processes. In fact, the robust process basically coincides with the nominal process in the nominal scenario for an uncertainty size of up to 5\%. Only at about 10\% uncertainty the objective reaches larger values than the nominal process in the nominal scenario. Even for the largest uncertainty of 20\%, the robust objective for the worst robust scenario is at most about a factor of two worse than that of the nominal process in the nominal setting. This shows that the cost of robustness, i.e., the loss of protection against uncertainty, is low. 

\begin{figure}[h]
\centering
\includegraphics[width=0.49\textwidth]{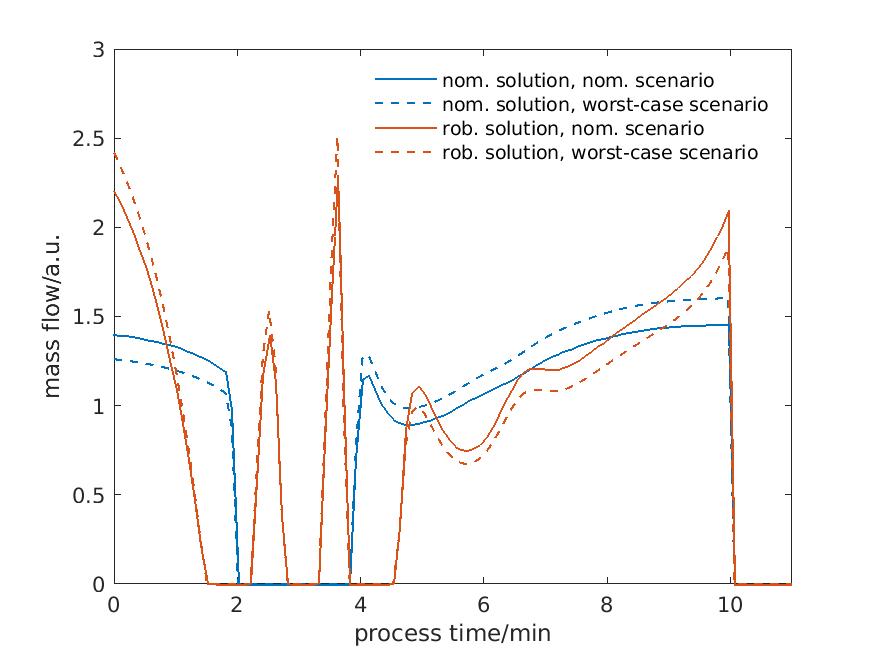}
\caption{The mass flow for the nominal process (blue), for the nominal scenario and the nominal worst-case scenario, and for the robust synthesis process (red), for the nominal scenario and the nominal worst-case scenario. The nominal worst-case scenario is a change of  the factor from $0.9$ to $1.1$ at time $2.2$ and the robust  worst-case scenario is a change of  the factor from $1.1$ to $0.9$ at time $3.9$.}
\label{fig:rob3}
\end{figure} 

Having considered the values of the objective function, we now take a closer look at the processes over time. Figure \ref{fig:rob3} shows the mass flows as a function of time in the optimized nominal and robust processes in their respective nominal and worst-case scenarios. It is clear that the nominal (and robust) solutions show similar patterns in both the nominal (robust) and worst-case scenarios. However, there are notable differences when comparing the nominal and robust processes. This is most prominent during process times between two and four minutes, where the worst-case scenarios also manifest themselves (nominal worst-case  at time 2.2, robust worst case at time 3.9). In particular, the robust process shows two additional peaks around these times that are not present in the nominal process.  Consequently, the robust process displays distinct characteristics that set it apart from the nominal, unprotected process.   

\begin{figure}[h]
\centering
\includegraphics[width=0.49\textwidth]{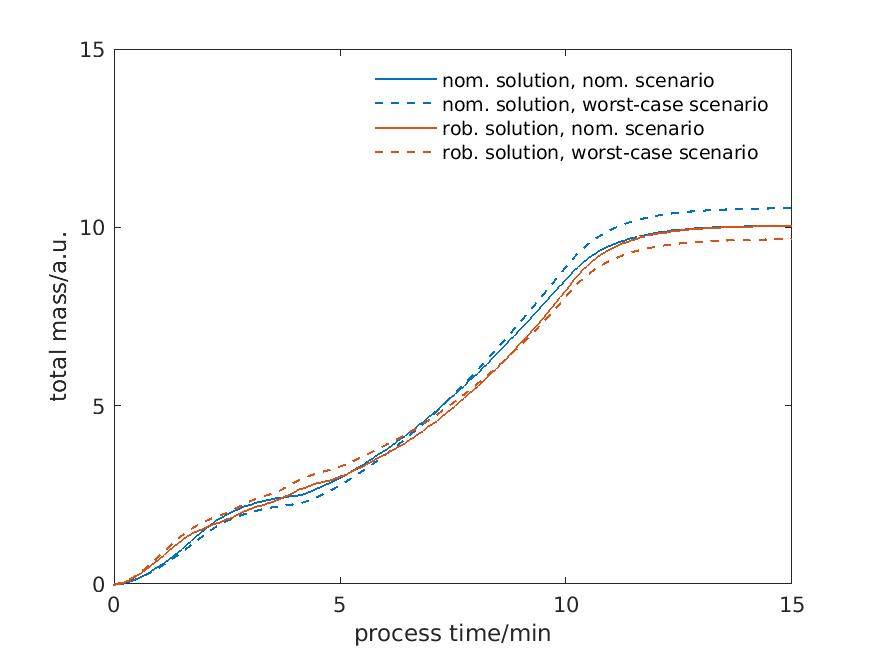}
\caption{The total mass for the nominal process (blue), for the nominal scenario and the nominal worst-case scenario, and for the robust process (red), for the nominal scenario and the nominal worst-case scenario. }
\label{fig:rob4}
\end{figure} 

Figure \ref{fig:rob4} illustrates the increase in total mass over time during the synthesis process. It can be observed that the total mass increases continuously until it reaches a saturation point at about 10 minutes, after which it remains relatively constant.

\begin{figure}[h]
\centering
\includegraphics[width=0.49\textwidth]{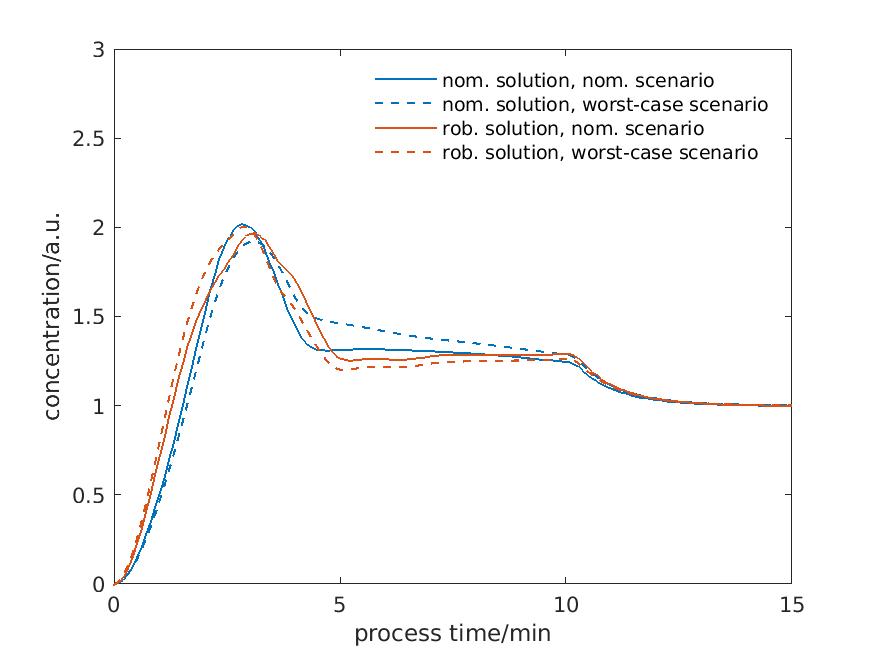}
\caption{The concentration for the nominal process (blue), for the nominal scenario and the nominal worst-case scenario, and for the robust process (red), for the nominal scenario and the nominal worst-case scenario. }
\label{fig:rob5}
\end{figure} 

In Figure \ref{fig:rob5} we present the concentration dynamics of the product over time during the process. Initially, the concentration rises rapidly and reaches a peak at approximately three minutes.  It then decreases and stabilises at a plateau between five and ten minutes. After ten minutes, the concentration continues to fall and asymptotically reach the equilibrium concentration of 1.

\begin{figure}[h]
\centering
\includegraphics[width=0.49\textwidth]{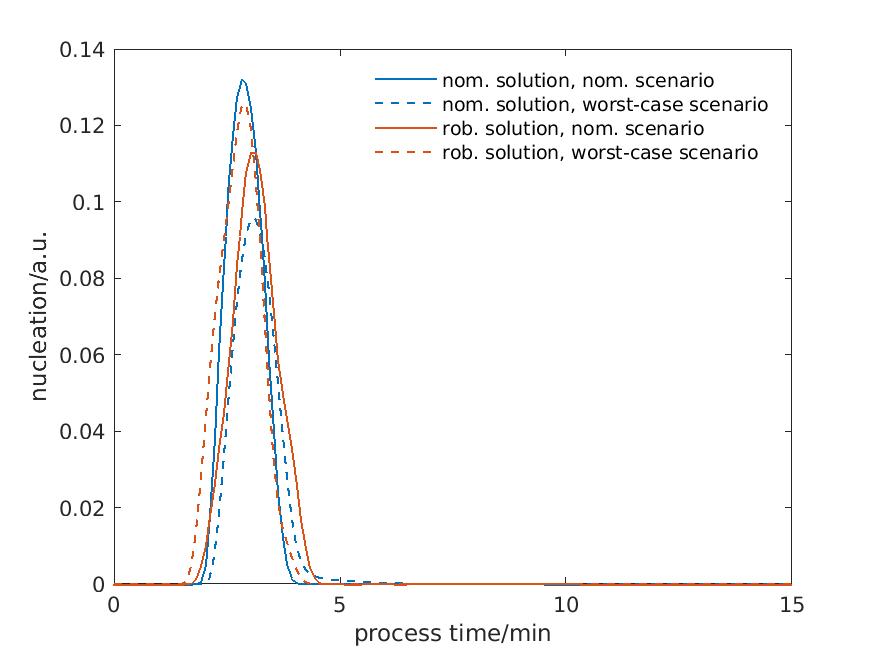}
\caption{The nucleation for the nominal process (blue), for the nominal scenario and the nominal worst-case scenario, and for the robust process (red), for the nominal scenario and the nominal worst-case scenario. }
\label{fig:rob6}
\end{figure} 
To study nucleation and growth, the nucleation function is shown in Figure \ref{fig:rob6}. Whereas it takes some time before nucleation takes place, nucleation is present in the range of 2 to 5 minutes, exhibits a sharp peak around 3 minutes, and remains approximately at zero after about 5 minutes.  The robust solution starts nucleating a bit later and less pronounced, when compared to the nominal process. 

\begin{figure}[h]
\centering
\includegraphics[width=0.49\textwidth]{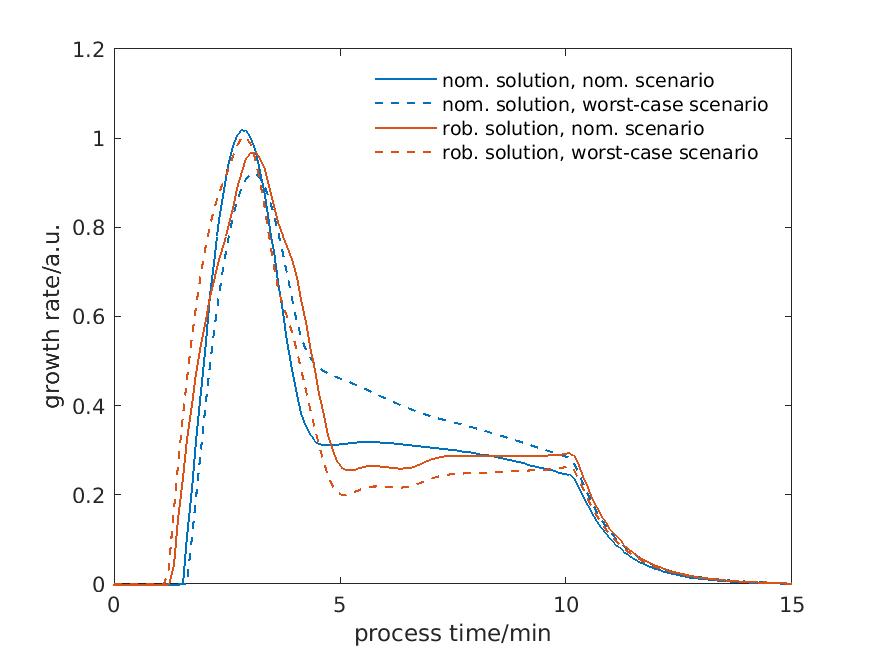}
\caption{The growth rate for the nominal process (blue), for the nominal scenario and the nominal worst-case scenario, and for the robust process (red), for the nominal scenario and the nominal worst-case scenario. }
\label{fig:rob7}
\end{figure} 

To study the growth of the particle after nucleation has taken place, Figure \ref{fig:rob7} shows the growth rate over time. It initially increases rapidly and reaches a peak at about 3 minutes. It then declines until it reaches roughly a plateau at about 5 minutes. Thereafter, the particles continue to grow at a roughly constant rate until about 10 minutes of process time. The growth rates differ slightly between the nominal and the robust process between 5 and 10 minutes. As the nominal growth rate slowly decreases slightly, the robust growth rate increases slightly. Subsequently, all growth rates decrease rapidly to zero,  resulting in the end product having reached its final size distribution.

\begin{figure}[h]
\centering
\includegraphics[width=0.49\textwidth]{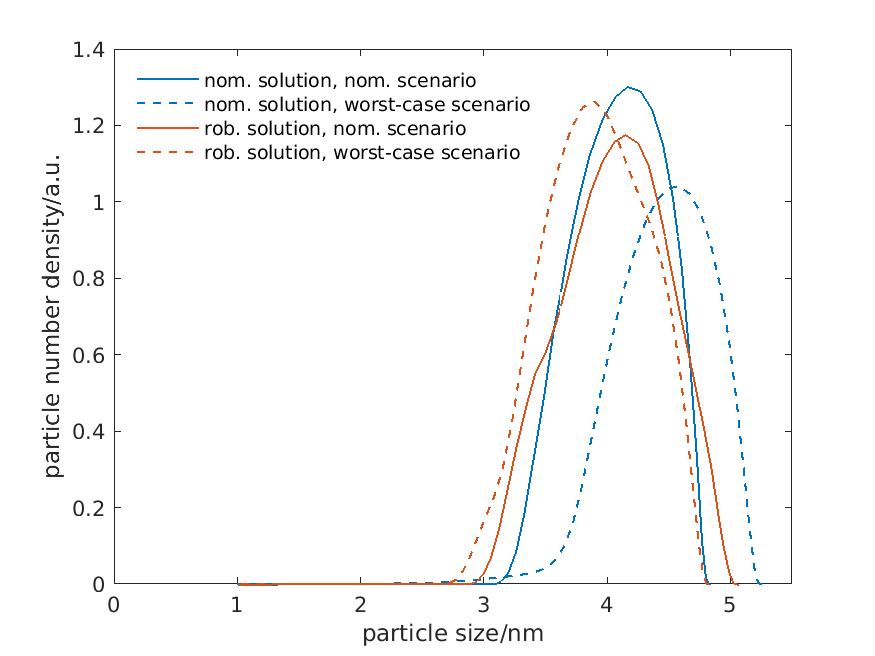}
\caption{The particle size distribution for the nominal process (blue), for the nominal scenario and the nominal worst-case scenario, and for the robust process (red), for the nominal scenario and the nominal worst-case scenario. }
\label{fig:rob8}
\end{figure} 

Finally, we examine the particle size distributions in Figure \ref{fig:rob8}, which shows the particle number density as a function of particle size. The objective of our optimization problem is to achieve a mean particle size of 4 nm, which is successfully achieved by both the nominal and the robust process in the nominal scenario. However, this changes drastically when evaluating the processes in disturbed settings. The robust process then produces nanoparticles with an average size of 3.85 nm, even in the worst case of uncertainty. This mean value is therefore still close to the desired size of 4 nm, even under uncertainty. In comparison, applying the nominal process in the nominal worst case results in significantly larger particles (4.58 nm). This yields an end product that is larger than the desired mean, making it potentially useless. In contrast, the robust process achieves the desired mean size. This clearly demonstrates the superiority of the robust process, which protects against uncertainties, over the unprotected nominal process. 

Finally, we study the price that we need to pay when implementing the robust process, which is typically referred to as the price of robustness. The latter here means to compare the nominal process with the robust process when however applied in the unperturbed nominal scenario, i.e., in the ideal case of perfect laboratory conditions with no uncertainties. The mean values of the PSDs almost coincide and are both about 4.15 nm for the nominal and robust process. This means that with respect to the means of the PSDs, there is no cost of robustness. However, when looking at their variances, it can be seen that the robust process leads to a slightly larger variance than that of the nominal process when the latter is used in the nominal scenario. In summary, this demonstrates the low case of robust protection and its obvious advantages when using the robust solution.    

\section{Conclusions and Outlook} \label{Subsec:conclusion}
In this work, we have introduced a robust optimization framework for particle synthesis, including nucleation and growth, modeled with population balance equations. We have formulated an optimization problem that focuses on the time-dependent precursor inflow rate. It minimizes the quadratic deviation from a given mean particle size and its variance. The resulting robust minimax problem assumes that there are uncertainties in the inflow rate and hedges against them. The numerical optimization was enabled by an advanced integration of a recently introduced adaptive bundle framework, which can solve nonconvex and nonlinear robust optimization problems, together with the exact method of moments. The algorithm was applied to particle synthesis processes with realistic parameter values taken from quantum dot synthesis. The computational results showed that the bundle method delivers good results in a short time while keeping the price of robustness low.  

We have investigated the nominal and robust processes, specifically evaluating them in their worst-case scenarios. It was found that the nominal solution is highly sensitive to disturbances and is unable to synthesize a suitable product under uncertainties. In contrast, the robust process successfully produces particles of the desired size at a relatively low cost of robustness. This is particularly evident in Figure \ref{fig:rob2}, where the cost of the robust solution is only marginally worse than that of the nominal solution, whereas the cost of the nominal solution is very large in the nominal worst case. This is even more striking in Figure \ref{fig:rob8}, where it is clearly visible that the robust process leads to better PSDs, while simultaneusly keeping the price of robust protection low. In fact, even in the nominal scenario,  which represents unperturbed, ideal laboratory conditions, the robust procedure leads to similar mean PSDs compared to the nominal solution, while the robust variance is only slightly larger than that of the nominal one.   
The computational results show that robust optimization of particle synthesis under uncertainty is crucial for the production of high-quality particles. We emphasize that the algorithm presented here is able to compute robust optimal processes in a short time. 

In the future, the adaptive bundle framework will be further extended by adaptive discretization widths with the aim of allowing very fine discretizations or continuous uncertainty sets. In conclusion, our computational results show that the adaptive bundle method, together with the exact methods of moments, is an efficient tool for dealing with uncertainties in the nanoparticle synthesis process, which can be computed quickly and at a low cost in terms of robustness. Finally, hedging mathematical optimization problems against uncertainty, especially when partial differential equations (PDEs) are involved, is also very relevant also in other contexts. Applications are numerous, ranging from particle problems to energy problems to problems in the natural sciences. The bundle method integrated by PDE-related solution approaches, as presented here, is not limited to particle synthesis. Thus, the approach to particle precipitation presented here will be developed for other such challenging problems in the future.  

\textit{CRediT author statement}\\
Martina Kuchlbauer: Methodology, Software, Validation, Investigation, Writing, Visualization\\
Jana Dienstbier: Methodology, Writing - Review \& Editing\\
Adeel Muneer: Software, Writing - Review \& Editing\\
Hanna Hedges: Software, Writing - First Draft\\
 Michael Stingl: Conceptualization, Methodology, Investigation, Validation, Writing - Review \& Editing\\
Frauke Liers: Conceptualization, Methodology, Validation, Investigation, Writing, Supervision\\
Lukas Pflug: Conceptualization, Software, Methodology, Validation, Formal Analysis, Investigation, Writing, Visualization

\textit{Funding}\\
This paper has been funded by the Deutsche Forschungsgemeinschaft
(DFG, German Research Foundation) - \linebreak Project-ID 416229255 -
SFB 1411 and by Project B06 in CRC TRR 154. 

\bibliographystyle{abbrv}
\bibliography{literatur}

\section{Appendix}

\subsection{Additional Computational Results} \label{sec:appendix-computations}

\begin{figure}[h]
\centering
\includegraphics[width=0.49\textwidth]{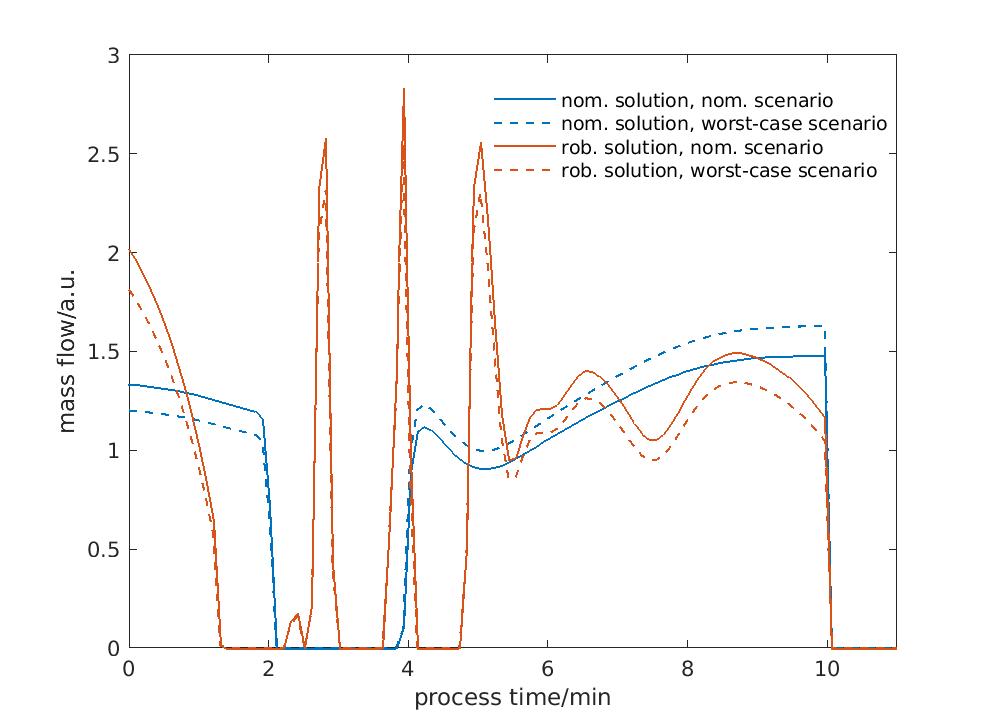}
\caption{For the case of $\omega_1 = 10$, $\omega_2=1$: The mass flow for the nominal solution (blue), for the nominal scenario and the nominal worst-case scenario, and for the robust solution (red), for the nominal scenario and the nominal worst-case scenario. }
\label{fig:rob1}
\end{figure} 

\begin{figure}[h]
\centering
\includegraphics[width=0.49\textwidth]{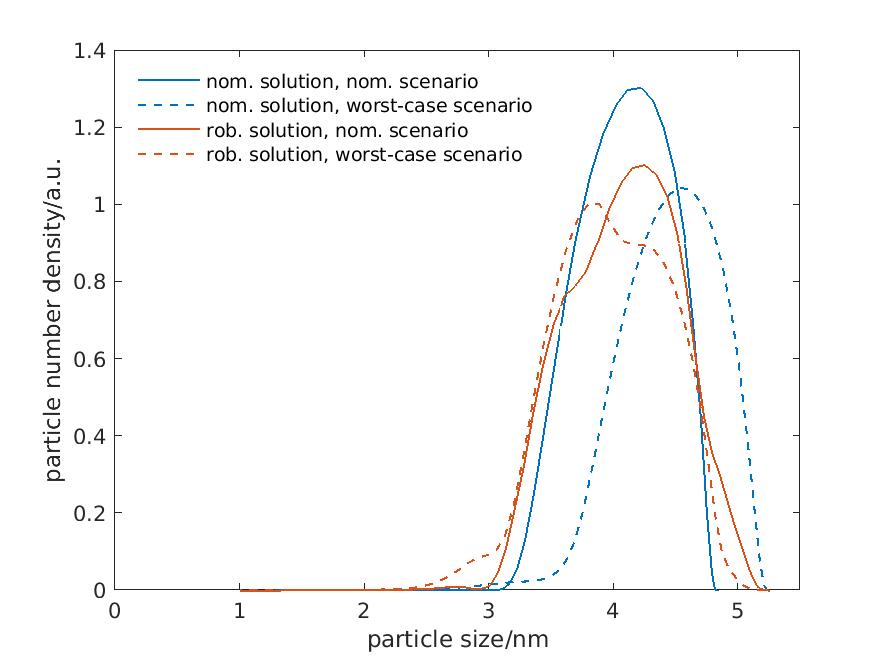}
\caption{For the case of $\omega_1 = 10$, $\omega_2=1$: The particle size distribution for the nominal solution (blue), for the nominal scenario and the nominal worst-case scenario, and for the robust solution (red), for the nominal scenario and the nominal worst-case scenario. }
\label{fig:rob9}
\end{figure}

\end{document}